\documentclass[11pt]{article}
\usepackage{amsfonts}
\usepackage{latexsym}
\usepackage{amsmath}
\usepackage{amssymb}
\usepackage{amsthm}
\newtheorem{theorem}{Theorem}[section]
\newtheorem{lemma}[theorem]{Lemma}

\newtheorem{proposition}[theorem]{Proposition}

\newtheorem{remark}[theorem]{Remark}
\theoremstyle{definition}
\newtheorem{definition}[theorem]{Definition}

\newtheorem*{notation}{Notation}

\theoremstyle{remark}
\newtheorem*{note*}{Note}

\makeatletter
\newcommand{\ls}{\leqslant}
\newcommand{\gr}{\geqslant}

\makeatother

\begin{document}
\small

\title{\bf A note on subgaussian estimates for linear
functionals on convex bodies}

\author{A. Giannopoulos\footnote{Research supported by the EPEAEK II program ``Pythagoras
II".} \quad A. Pajor\footnote{The second named author would like to
thank the Department of Mathematics of the University of Athens for
the hospitality.}\quad and\quad G. Paouris\footnote {The third named
author is supported by a Marie Curie Intra-European Fellowship
(EIF), Contract MEIF-CT-2005-025017.}}

\date{}

\maketitle

\begin{abstract}
\footnotesize We give an alternative proof of a recent result of
Klartag on the existence of almost subgaussian linear functionals on
convex bodies. If $K$ is a convex body in ${\mathbb R}^n$ with
volume one and center of mass at the origin, there exists $x\neq 0$
such that
$$|\{ y\in K:\,|\langle y,x\rangle |\gr t\|\langle\cdot
,x\rangle\|_1\}|\ls\exp (-ct^2/\log^2(t+1))$$ for all $t\gr 1$,
where $c>0$ is an absolute constant. The proof is based on the study
of the $L_q$--centroid bodies of $K$. Analogous results hold true
for general log-concave measures.
\end{abstract}

\section{Introduction}

The purpose of this note is to provide an alternative proof of a
recent result of Klartag (see \cite{Kl3}) on the existence of almost
subgaussian linear functionals on convex bodies. Let $K$ be a convex
body in ${\mathbb R}^n$ with volume $|K|=1$ and center of mass at
the origin. Let $\psi :[0,\infty )\to [0,\infty )$ be a convex,
increasing function with $\psi (0)=0$. For every bounded measurable
function $f:K\rightarrow {\mathbb R}$, define
$$\| f\|_{\psi }=\inf\left\{ t>0:\int_K
\psi (|f(x)|/t)\,dx\ls 1\right\}.\leqno (1.1)$$ We will be
interested in the $\psi_{\alpha }$--norm of linear functionals
$y\mapsto \langle y,x\rangle $ on $K$, where $1\ls\alpha\ls 2$ and
$\psi_{\alpha }(t)=e^{t^{\alpha }}-1$. We say that $x\neq 0$ defines
a $\psi_{\alpha }$--direction for $K$ with constant $B>0$ if
$$\|\langle\cdot ,x\rangle\|_{\psi_{\alpha }}\ls
B\|\langle \cdot ,x\rangle\|_1. \leqno (1.2)$$ It is not hard to
check that this holds true if and only if
$$\| \langle\cdot ,x\rangle\|_q\ls
cBq^{1/\alpha }\| \langle\cdot ,x\rangle\|_1\leqno (1.3)$$ for
every $q\gr 1$, where $c>0$ is an absolute constant. By Borell's
lemma (see \cite{MS}, Appendix III), there exists an absolute
constant $C>0$ such that if $K$ is a convex body in ${\mathbb
R}^n$, then every $x\neq 0$ is a $\psi_1$--direction for $K$ with
constant $C$.

The study of $\psi_2$--directions for linear functionals on convex
bodies is motivated by the study of isotropic convex bodies and
Bourgain's approach to the isotropic constant problem. A convex body
$K$ in ${\mathbb R}^n$ is called isotropic if it has volume $|K|=1$,
center of mass at the origin, and there exists a constant $L_K>0$
such that
$$\int_K\langle y,\theta\rangle^2dy =L_K^2\leqno (1.4)$$
for every $\theta\in S^{n-1}$. Every convex body with center of mass
at the origin has a linear image which is isotropic (see \cite{MP}).
This image is unique up to orthogonal transformations, and hence,
the isotropic constant $L_K$ is well--defined for the linear class
of $K$. The isotropic constant problem asks if there exists an
absolute constant $C>0$ such that $L_K\ls C$ for every isotropic
convex body in any dimension. One can easily see that
$L_K=O(\sqrt{n})$ for every $K$. Uniform boundedness of $L_K$ is
known for some classes of bodies: unit balls of spaces with
$1$--unconditional basis, zonoids and their polars, etc. Bourgain
(see \cite{Bou1}) proved that $L_K=O(\sqrt[4]{n}\log n)$ and, very
recently, Klartag (see \cite{Kl2}) improved this bound to
$L_K=O(\sqrt[4]{n})$. Moreover, in \cite{Bou2} Bourgain proved that
if every $x\neq 0$ is a $\psi_2$--direction for $K$ with constant
$B$, then $L_K$ is bounded by $cB\log (B+1)$.

A question of Milman, related to this line of thought, is whether,
for every isotropic convex body $K$ in ${\mathbb R}^n$, most
$\theta\in S^{n-1}$ define a $\psi_2$--direction for $K$ with a
``good" constant (for example, logarithmic in $n$). Until recently,
it was not known if there exists an absolute constant $C>0$ such
that every isotropic convex body has at least one
$\psi_2$--direction with constant $C$. Some positive results are
known for special classes of convex bodies. Bobkov and Nazarov (see
\cite{BN1} and \cite{BN2}) have proved that if $K$ is an isotropic
$1$--unconditional convex body, then $\|\langle\cdot
,x\rangle\|_{\psi_2} \ls c\sqrt{n}\| x\|_{\infty }$ for every $x\neq
0$. This shows that the diagonal direction is a $\psi_2$--direction.
For the class of zonoids, the existence of good $\psi_2$--directions
was established in \cite{Pa1}. Another partial result, which gives
more information in the case of isotropic convex bodies with ``small
diameter", was obtained in \cite{Pa2}: If $K\subseteq
(\gamma\sqrt{n}L_K)B_2^n$ for some $\gamma >0$, then
$$\sigma \big (\theta\in S^{n-1}:\|\langle \cdot ,\theta
\rangle\|_{\psi_2}\gr c_1\gamma tL_K\big )\ls\exp
(-c_2\sqrt{n}t^2/\gamma )\leqno (1.5)$$ for every $t\gr 1$, where
$\sigma $ is the rotationally invariant probability measure on
$S^{n-1}$ and $c_1,c_2>0$ are absolute constants.

Klartag (see \cite{Kl3}) gave a positive answer to this question,
showing that every isotropic convex body admits at least one almost
subgaussian linear functional. Our aim is to give a second (short)
proof of this fact.

\begin{theorem}
Let $K$ be an isotropic convex body in ${\mathbb R}^n$. There exists
$x\neq 0$ such that
$$|\{ y\in K:\,|\langle y,x\rangle |\gr t\|\langle\cdot
,x\rangle\|_1\}|\ls\exp (-ct^2/\log^{\tau }(t+1))\leqno (1.6)$$
for all $t\gr 1$, where $c,\tau >0$ are absolute constants.
\end{theorem}

It is clear that if $x$ defines a $\psi_{\alpha }$--direction for
$K$ and if $T\in SL(n)$, then $T^{\ast }x$ defines a $\psi_{\alpha
}$--direction (with the same constant) for $T(K)$. It follows that
Theorem 1.1 provides almost subgaussian directions for every convex
body: If $K$ is a convex body in ${\mathbb R}^n$ with volume one and
center of mass at the origin, there exists $x\neq 0$ such that (1.6)
holds true for all $t\gr 1$.

The argument of Klartag is based on the study of the level sets of
the logarithmic Laplace transform of log--concave functions. The
argument we present here is based on the study of the
$L_q$--centroid bodies of an isotropic convex body. This family of
bodies was studied and used by the third named author in
\cite{Pa2}, and in particular in \cite{Pa3}, where the following
sharp dimension--dependent concentration of volume estimate was
proved: There exists an absolute constant $c>0$ such that if $K$
is an isotropic convex body in ${\mathbb R}^n$, then
$$\left |\big\{ x\in K:\| x\|_2\gr  c\sqrt{n}L_Kt\big \}\right |
\ls \exp\left ( -\sqrt{n}t\right )\leqno (1.7)$$ for every $t\gr
1$, where $\|\cdot\|_2$ is the Euclidean norm. The tools which are
developed in \cite{Pa3} allow us to give a very simple proof of
Theorem 1.1. We present an argument which gives $\tau =2$, i.e.
the upper bound in (1.6) is $\exp (-ct^2/\log^{2}(t+1))$.

\medskip

\begin{notation}We work in ${\mathbb R}^n$, which is
equipped with a Euclidean structure $\langle\cdot ,\cdot\rangle $.
We denote by $\|\cdot \|_2$ the corresponding Euclidean norm, and
write $B_2^n$ for the Euclidean unit ball, and $S^{n-1}$ for the
unit sphere. Volume is denoted by $|\cdot |$. If $K$ is a convex
body in ${\mathbb R}^n$, we set $\overline{K}=K/|K|^{1/n}$; this is
the dilation of $K$ which has volume one. We write $\sigma $ for the
rotationally invariant probability measure on $S^{n-1}$. The
Grassmann manifold $G_{n,k}$ of $k$--dimensional subspaces of
${\mathbb R}^n$ is equipped with the Haar probability measure
$\mu_{n,k}$.

A convex body is a compact convex subset $C$ of ${\mathbb R}^n$
with non--empty interior. We say that $C$ has center of mass at
the origin if $\int_C\langle x,\theta\rangle dx=0$ for every
$\theta\in S^{n-1}$. The support function $h_C:{\mathbb
R}^n\rightarrow {\mathbb R}$ of $C$ is defined by $h_C(x
)=\max\{\langle x,y\rangle :y\in C\}$. The mean width of $C$ is
defined by $$w(C)=\int_{S^{n-1}}h_C(\theta )\,\sigma (d\theta
).\leqno (1.8)$$ The letters $c,c^{\prime }, c_1, c_2$ etc. denote
absolute positive constants which may change from line to line. We
refer to the books \cite{Schn}, \cite{MS} and \cite{Pi} for basic
facts from the Brunn--Minkowski theory and the asymptotic theory
of finite dimensional normed spaces.
\end{notation}

\section{Normalized $L_q$--centroid bodies}

Let $K$ be a convex body of volume $1$ in ${\mathbb R}^n$. For every
$q\gr 1$ we define the $L_q$--centroid body $Z_q(K)$ of $K$ by its
support function:
$$h_{Z_q(K)}(x)=\|\langle \cdot ,x\rangle\|_q:=
\left (\int_K|\langle y,x\rangle |^qdy\right )^{1/q}. \leqno (2.1)$$
Since $|K|=1$, we readily see that $Z_1(K)\subseteq Z_p(K)\subseteq
Z_q(K)\subseteq Z_{\infty }(K)$ for every $1\ls p\ls q\ls \infty $,
where $Z_{\infty }(K)={\rm conv}\{ K,-K\}$. On the other hand, one
has the reverse inclusions
$$Z_q(K)\subseteq \frac{cq}{p}Z_p(K)\leqno (2.2)$$
for every $1\ls p< q<\infty $, as a consequence of the
$\psi_1$--behavior of $y\mapsto\langle y,x\rangle $. Observe that
$Z_q(K)$ is always symmetric, and $Z_q(TK)=T(Z_q(K))$ for every
$T\in SL(n)$ and $q\in [1,\infty ]$. Also, if $K$ has its center of
mass at the origin, then $Z_q(K)\supseteq cZ_{\infty }(K)$ for all
$q\gr  n$, where $c>0$ is an absolute constant.

It should be mentioned that $L_q$--centroid bodies were introduced
in \cite{LZ} under a different normalization. Lutwak, Yang and Zhang
(see \cite{LYZ} and \cite{CG} for a different proof) have
established the $L_q$ affine isoperimetric inequality
$$|Z_q(K)|^{1/n}\gr |Z_q(\overline{B_2^n})|^{1/n} \gr c\sqrt{q/n}\leqno (2.3)$$
for every $1\ls q\ls n$, where $c>0$ is an absolute constant.

We will need upper estimates for the quermassintegrals of the
$L_q$--centroid bodies of an isotropic convex body. These follow
immediately from estimates on the projections of $Z_q(K)$, which are
obtained in \cite{Pa3}. Fix $1\ls k\ls n$ and a $k$--dimensional
subspace $F$ of $\mathbb R^n$, and denote by $E$ the orthogonal
subspace of $F$. For every $\phi \in S_{F}$, define $E(\phi)=\{y\in
{\rm span}\{E, \phi \}: \langle y,\phi\rangle \gr 0 \}$. By a
theorem of K. Ball (see \cite{Ball} and \cite{MP}), for every convex
body $K$ of volume one in $\mathbb R^n$, for every $q\gr  0$ and
every $\phi \in F$, the function
$$\phi\mapsto \|\phi\|_2^{1+\frac{q}{q+1}}\left( \int_{K \cap E(\phi)} |\langle y, \phi \rangle |^q
dy
\right)^{-\frac{1}{q+1}}\leqno (2.4)$$ is a gauge function on $F$
(see also \cite{BKM} for the not necessarily symmetric case). If we
denote by $B_q(K,F)$ the convex body in $F$ whose gauge function is
defined by (2.4), then the volume of $B_q(K,F)$ is given by
$$|B_q (K,F) | = |B_2^k| \int_{S_F} \left( \int_{K \cap E(\phi)}
 |\langle x, \phi \rangle |^{q}  dx \right)^{\frac{k}{q+1}} d\sigma_F (\phi).\leqno (2.5)$$
The following identity was proved in \cite{Pa3}.

\begin{proposition}
Let $K$ be a convex body of volume $1$ in ${\mathbb
R}^n$ and let $1\ls  k\ls  n-1$. For every $F\in G_{n,k}$ and
every $q\gr 1$ we have that
$$P_F (Z_q(K)) = \left (k+q\right )^{1/q} |B_{k+q-1} (K,F)|^{1/k +
1/q} Z_q (\overline{B}_{k+q-1} (K,F)).\leqno (2.6)$$
\end{proposition}

Using this identity and exploiting (2.5) in order to estimate the
volume of $B_q(K,F)$, one gets the following estimate (see
\cite{Pa3}).

\begin{proposition}
Let $K$ be an isotropic convex body in ${\mathbb R}^n$. If $F\in
G_{n,k}$ and $E=F^{\perp }$ then, for every $q\in {\mathbb N}$ we
have that
$$P_{F}(Z_q(K)) \subseteq \frac{c(k+q)}{k} L_K
Z_q(\overline{B}_{k+q-1}(K,F))\leqno (2.7)$$ where $c>0$ is an
absolute constant.
\end{proposition}

\begin{definition}{\rm Let $K$ be an isotropic convex body in ${\mathbb R}^n$.
For every integer $q\gr 1$ we define the {\it normalized
$L_q$--centroid body} $K_q$ of $K$ by}
$$K_q=\frac{1}{\sqrt{q}L_K}Z_q(K).\leqno (2.8)$$
\end{definition}

Since $|Z_q(\overline{B}_{k+q-1}(K,F))|\ls
|\overline{B}_{k+q-1}(K,F)|=1$, Proposition 2.2 shows that
$$|P_F(K_q)|^{1/k}\ls
\frac{c(k+q)}{k\sqrt{q}}|Z_q(\overline{B}_{k+q-1}(K,F))|^{1/k}\ls
\frac{c_1(k+q)}{k}\frac{\sqrt{k}}{\sqrt{q}}|B_2^k|^{1/k}\leqno
(2.9)$$ for every $F\in G_{n,k}$. If $1\ls k\ls q$, this estimate
takes the simpler form
$$|P_F(K_q)|^{1/k}\ls 2c_1\frac{\sqrt{q}}{\sqrt{k}}|B_2^k|^{1/k}.\leqno
(2.10)$$ In particular, for every $F\in G_{n,q}$ we have
$$|P_F(K_q)|^{1/k}\ls 2c_1|B_2^k|^{1/k}.\leqno
(2.11)$$ A standard argument (based on the log--concavity of the
quermassintegrals of $P_F(K_q)$) implies that since (2.11) is true
for every $F\in G_{n,q}$, it remains valid for every $F\in G_{n,k}$,
where $q\ls k\ls n$. We summarize these observations in the next
Theorem.

\begin{theorem}
Let $K$ be an isotropic convex body in ${\mathbb R}^n$. If $1\ls
k,q\ls n$ are integers, and if $F\in G_{n,k}$, then
$$|P_F(K_q)|^{1/k}\ls c_1\max \{ \sqrt{q/k}, 1\}|B_2^k|^{1/k},\leqno (2.12)$$
where $c_1>0$ is an
absolute constant. In particular,
$$|K_q|^{1/n}\ls  c_1|B_2^n|^{1/n}.\leqno (2.13)$$
\end{theorem}

The last ingredient of the proof is a consequence of the main
result in \cite{Pa3}: from (1.7) it follows that
$$\left (\int_K\| y\|_2^q\,dy\right )^{1/q}\ls c\sqrt{n}L_K\leqno
(2.14)$$ for all $1\ls q\ls\sqrt{n}$. Since
$$w(Z_q(K))\leq\left (\int_{S^{n-1}}\int_K|\langle y,\theta\rangle
|^qdy\,\sigma (d\theta )\right )^{1/q}\ls \left (
\frac{C\sqrt{q}}{\sqrt{n}}\int_K\| y\|_2^q\,dy\right )^{1/q}\leqno
(2.15)$$ for all $1\ls q\ls n$, we have the following Lemma.

\begin{lemma}
Let $K$ be an isotropic convex body in ${\mathbb R}^n$. If $1\ls
q\ls\sqrt{n}$, then
$$w(K_q)\ls C,\leqno (2.16)$$ where $C>0$ is an absolute constant.
\end{lemma}

\begin{remark}{\rm Without using Lemma 2.5, which fully exploits
the results of \cite{Pa3}, we can prove Theorem 1.1 with $\tau
=2+\epsilon $ for any $\epsilon >0$.}
\end{remark}

\section{Covering numbers of $K_q$}

Let $N(K_q,sB_2^n)$ denote the minimal number of translates of
$sB_2^n$ whose union covers $K_q$. A standard way to estimate the
covering number $N(K_q,sB_2^n)$ is through the inequality
$$|tB_2^n|\cdot N(K_q,2tB_2^n)\ls |K_q+tB_2^n|,\leqno (3.1)$$ which is
valid for every $t>0$. We will use our information on the
projections of $K_q$ in order to give an upper bound for
$|K_q+tB_2^n|$.

\begin{proposition}Let $K$ be an isotropic convex body in ${\mathbb
R}^n$. For every $1\ls q\ls n$ and every $t>0$, we have that
$$N(K_q,2tB_2^n)\ls \exp \left ( C\frac{\sqrt{qn}}{\sqrt{t}}
+C\frac{n}{t}\right ),\leqno (3.2)$$ where $C>0$ is an absolute
constant.
\end{proposition}

\noindent {\it Proof}. From the classical Steiner's formula we know
that
$$|K_q+tB_2^n|=\sum_{k=0}^n{n\choose k}W_{[n-k]}(K_q)t^{n-k}\leqno (3.3)$$
for all $t>0$, where $W_{[n-k]}(K_q)$ is the mixed volume
$V_{k}(K_q)=V(K_q;k,B_2^n;n-k)$ (see \cite{Schn}).

We will use Kubota's integral formula to express $W_{[n-k]}(K_q)$ as
an average of the volumes of the $k$--dimensional projections of
$K_q$: for every $1\ls k\ls n-1$ we have
$$W_{[n-k]}(K_q) = \frac{|B_2^n|}{|B_2^k|} \int_{G_{n,k}} |P_F
(K_q)|\, d\mu_{n,k}(F).\leqno (3.4)$$ Using (3.3), (3.4) and the
estimates from Theorem 2.4, we can write
$$|K_q+tB_2^n|\ls |B_2^n|\sum_{k=0}^n{n\choose k}\left
(c_1\max \{ \sqrt{q/k}, 1\}\right )^kt^{n-k}.\leqno (3.5)$$ Then,
(3.1) shows that
$$N(K_q,2tB_2^n)\ls\sum_{k=0}^q\left ( \frac{c_2n\sqrt{q}}{k^{3/2}t}\right
)^k+\sum_{k=q+1}^n\left ( \frac{c_2n}{kt}\right )^k.\leqno (3.6)$$
Observe that for $1\ls k\ls q$ we have
$$\left ( \frac{c_2n\sqrt{q}}{k^{3/2}t}\right
)^k\ls \left (\frac{c_2nq}{k^2t}\right )^k\ls \frac{\left
(c_3\sqrt{nq/t}\right )^{2k}}{(2k)!},\leqno (3.7)$$ while, for
$q\ls k\ls n$ we have
$$\left ( \frac{c_2n}{kt}\right )^k\ls \frac{\left (c_4n/t\right )^k}{k!}.\leqno (3.8)$$
It follows that
$$N(K_q,2tB_2^n)\ls \exp \left (c_3\frac{\sqrt{qn}}{\sqrt{t}}\right )
+\exp\left ( c_4\frac{n}{t}\right ),\leqno (3.9)$$ and the result
follows. $\hfill\Box $

\begin{remark}{\rm The proof actually gives $N(K_q,2tB_2^n)\ls
\exp\left ( C\frac{n^{2/3}q^{1/3}}{t^{2/3}}+C\frac{n}{t}\right )$
for every $t>0$, but this would play no role in the proof of the
main result.}
\end{remark}

\section{Proof of the Theorem}

Let $K$ be an isotropic convex body in ${\mathbb R}^n$. Consider
the convex body
$$T={\rm conv}\left (\bigcup_{i=1}^{\lfloor\log_2n\rfloor
}\frac{1}{i}K_{2^i}\right ).\leqno (4.1)$$ We will use the
following standard fact.

\begin{lemma}Let $A_1,\ldots ,A_s$ be subsets of $RB_2^n$. For every $t>0$ we have
that
$$N({\rm conv}(A_1\cup\cdots\cup A_s),2tB_2^n)\ls \left
(\frac{cR}{t}\right )^s\prod_{i=1}^sN(A_i,tB_2^n).\leqno (4.2)$$
\end{lemma}

\noindent {\it Sketch of the proof.} For $i=1,\ldots ,s$, let
$N_i$ be a subset of ${\mathbb R}^n$ with cardinality
$|N_i|=N(A_i,tB_2^n)$, so that $A_i\subseteq\bigcup_{x_i\in
N_i}(x_i+tB_2^n)$. Let $B_1^s$ denote the unit ball of $\ell_1^s$
and fix $Z\subseteq B_1^s$ of minimal cardinality, so that
$B_1^s\subseteq \bigcup_{z\in Z}(z+(t/R)B_1^n)$. It is well--known
that $|Z|\ls (cR/t)^s$, where $c>0$ is an absolute constant.
Consider the set $N=\{ w=z_1x_1+\cdots +z_sx_s:\;x_i\in N_i,
z=(z_1,\ldots ,z_s)\in Z\}$. Then, ${\rm conv}(A_1\cup\cdots\cup
A_s)\subseteq\bigcup_{w\in N}(w+2tB_2^n)$. $\hfill\Box $

\medskip

Let $s=\lfloor\log_2n\rfloor $ and $m=\lfloor
\log_2(\sqrt{n})\rfloor \simeq s/2$. We apply Lemma 4.1 with
$A_i=\frac{1}{i}K_{2^i}$, $1\ls i\ls s$, and $t=1$. Observe that
$A_i\subseteq c_1\sqrt{n}B_2^n$ for all $i\ls s$ (to see this,
recall the known fact that if $K$ is an isotropic convex body in
${\mathbb R}^n$, then $K\subseteq (cnL_K)B_2^n$). Using Sudakov's
inequality (see \cite{Pi}) and Lemma 2.5 to estimate $N(A_i,B_2^n)$
for $i\ls m$, and using the entropy estimates of Section 3 to
estimate $N(A_i,B_2^n)$ for $m<i\ls s=\lfloor\log_2n\rfloor $, we
may write
\begin{eqnarray*} N(T,B_2^n) &\ls & (c_2\sqrt{n})^{\lfloor\log_2n\rfloor }
\left [\prod_{i=1}^{\lfloor\log_2n\rfloor }N(K_{2^i},iB_2^n)\right ]\\
&\ls & e^{c_3n}\exp \left (
C\sqrt{n}\sum_{i=s+1}^{\lfloor\log_2n\rfloor }2^{i/2}\right
)\times \exp \left ( Cn\cdot \left (\sum_{i=1}^{m}\frac{1}{i^{2}}+
\sum_{i=m+1}^{2m}\frac{1}{i}\right )\right )\\
&\ls & e^{cn}.
\end{eqnarray*}
It follows that $|T|\ls |CB_2^n|$, where $C>0$ is an absolute
constant. Therefore, there exists $x\neq 0$ such that
$$h_T(x)\ls C\| x\|_2,\leqno (4.3)$$
and hence,
$$\|\langle\cdot ,x\rangle\|_{2^i}\ls C\,2^{i/2}iL_K\| x\|_2\leqno (4.4)$$ for every $i=1,2,\ldots
,\lfloor\log_2n\rfloor $. This easily implies the following.

\begin{theorem}Let $K$ be an isotropic convex body in ${\mathbb
R}^n$. There exists $\theta\in S^{n-1}$ such that
$$\|\langle\cdot ,\theta\rangle\|_q\ls C\sqrt{q}\,\log q\|\langle
\cdot ,\theta\rangle\|_2\leqno (4.5)$$ for every $q\gr 2$, where
$C>0$ is an absolute constant.
\end{theorem}

A standard argument shows that Theorem 4.2 implies Theorem 1.1 (it
is actually equivalent to Theorem 1.1 with $\tau =2$).

\begin{remark}{\rm The proof of Theorem 4.2 carries over to the case of
an arbitrary log-concave measure: the approach of \cite{Pa3} and
all the arguments we have used in this note depend only on the
Brunn--Minkowski theory. It follows that if $\mu $ is an isotropic
log--concave measure in ${\mathbb R}^n$, then there exists
$\theta\in S^{n-1}$ such that
$$\|\langle\cdot ,\theta\rangle\|_{L^q(\mu )}\ls C\sqrt{q}\,\log q\|\langle
\cdot ,\theta\rangle\|_{L^2(\mu )}\leqno (4.6)$$ for all $2\ls
q\ls n$, where $C>0$ is an absolute constant.}
\end{remark}

\bigskip

\footnotesize
\bibliographystyle{amsplain}

\bigskip

\bigskip

\noindent \textsc{A.\ Giannopoulos}: Department of Mathematics,
University of Athens, Panepistimiopolis 157 84, Athens, Greece.

\smallskip

\noindent \textit{E-mail:} \texttt{apgiannop@math.uoa.gr}

\medskip

\noindent \textsc{A.\ Pajor}: \'Equipe d'Analyse et de
Math\'ematiques Appliqu\'ees, Universit\'e de Marne-la-Vall\'ee,
Champs sur Marne, 77454, Marne-la-Vall\'ee, Cedex 2, France.

\smallskip

\noindent \textit{E-mail:} \texttt{Alain.Pajor@univ-mlv.fr}

\medskip

\noindent \textsc{G.\ Paouris}: \'Equipe d'Analyse et de
Math\'ematiques Appliqu\'ees, Universit\'e de Marne-la-Vall\'ee,
Champs sur Marne, 77454, Marne-la-Vall\'ee, Cedex 2, France.

\smallskip

\noindent \textit{E-mail:} \texttt{grigoris\_paouris@yahoo.co.uk}

\end{document}